\documentclass[12pt]{amsart}

\usepackage{amsfonts,amsmath,amssymb,enumerate,
hyperref,latexsym,xcolor,amsthm,graphicx,multirow, threeparttable}
\def\a{\alpha}
\def\b{\beta}
\def\l{\langle} \def\r{\rangle}
\def\div{\,\,\big|\,\,}
\newcommand\Ga{\Gamma}

\def\mod{{\sf mod~}}
\newcommand\bv{\mathbf{v}}
\newcommand\bu{\mathbf{u}}\newcommand\bw{\mathbf{w}}
\newcommand\bx{\mathbf{x}}
\newcommand\bo{\mathbf{o}}\newcommand\ba{\mathbf{a}}
\newcommand\bd{\mathbf{d}}\newcommand\be{\mathbf{e}}
\newcommand\bb{\mathbf{b}}   \newcommand\bc{\mathbf{c}}
 \newcommand\bt{\mathbf{t}}
 \newcommand\wt{\mathsf{wt}}\newcommand\supp{\mathsf{supp}}

\newcommand\ZZ{\mathbb{Z}} \newcommand\FF{\mathbb{F}} 
\newcommand\A{\mathrm{A}} \newcommand\Sy{\mathrm{S}}
\newcommand\Aut{\mathrm{Aut}}  
 \newcommand\Cos{\mathrm{Cos}} 
\newcommand\soc{\mathrm{soc}}
\newcommand\bC{\mathbf{C}}\newcommand\bD{\mathbf{D}}
\newcommand\bI{\mathbf{I}}\newcommand\bA{\mathbf{A}}
\newcommand\cC{\mathcal{C}}\newcommand\bP{\mathbf{P}}
\newcommand\N{\mathbf{N}} 
 
\newcommand\GL{\mathrm{GL}}  \newcommand\AGL{\mathrm{AGL}} \newcommand\PGL{\mathrm{PGL}}  \newcommand\PSL{\mathrm{PSL}}

 \newcommand\Q{\mathrm{Q}}  
 \newcommand\Ree{\mathrm{Ree}}

\newtheorem{theorem}{Theorem}[section]

\newtheorem{lemma}[theorem]{Lemma}

\theoremstyle{definition}
\newtheorem{example}[theorem]{Example}

\begin{document}
\title[2-arc-transitive graphs]{On $2$-arc-transitive graphs of product action type}
\thanks{Supported by the National Natural Science Foundation of China
(11971248), and the Fundamental Research Funds for the Central
Universities.}

\author[Z.P. Lu]{Zai Ping Lu}
\address{Z.P. Lu\\ Center for Combinatorics, LPMC\\ Nankai University\\ Tianjin 300071\\ P. R. China}
\email{lu@nankai.edu.cn}

\begin{abstract}
 In this paper, we discuss the structural information about $2$-arc-transitive (non-bipartite and bipartite) graphs of product action type. It is proved that a $2$-arc-transitive graph of product action type requires certain restrictions on either the vertex-stabilizers or the valency.  Based on the existence of  some equidistant linear codes, a construction is given for  $2$-arc-transitive graphs of non-diagonal product action type, which produces several families of such graphs. Besides, a nontrivial construction is given for  $2$-arc-transitive bipartite graphs of diagonal product action type
 %Besides, some examples are constructed.

\vskip 5pt

\end{abstract}
\maketitle
%%------------------------section----------------------------

%\parskip 5pt
\baselineskip 15pt

\vskip 20pt

\section{Introduction}\label{intro}
All graphs considered in this paper are assumed to be finite, simple
and undirected.

\vskip 5pt

Let $\Ga=(V,E)$ be a connected graph with vertex set $V$ and edge set $E$. An arc in $\Ga$ is an ordered pair of  adjacent vertices, and a $2$-arc is a triple $(\a,\b,\gamma)$ of distinct vertices with $\{\a,\b\},\{\b,\gamma\}\in E$. Denote by $\Aut(\Ga)$ the full automorphism group of $\Ga$. For a subgroup $G\leqslant \Aut(\Ga)$, the graph $\Ga$ is said to be $(G,2)$-arc-transitive   (or $(G,2)$-arc-regular)
if $G$ acts transitively (or regularly) on the set of $2$-arcs of $\Ga$, while the group $G$ is called
 a $2$-arc-transitive (or $2$-arc-regular) group of $\Ga$.
% (Note, the $2$-arc-transitivity of $G$ yields the
 %transitivity of $G$ on  the vertex set and  
 %the arc set of $\Ga$.)
For a vertex $\a\in V$, let $G_\a=\{g\in G\mid \a^g=\a\}$ and $\Ga(\a)=\{\b\in V\mid \{\a,\b\}\in E\}$, called the stabilizer   and the neighborhood of $\a$ in $G$ and $\Ga$, respectively. It is well-known that $G$ is $2$-arc-transitive if and only if $G$ acts transitively on $V$ and, for $\a\in V$, the stabilizer $G_\a$ acts $2$-transitively on $\Ga(\a)$.

\vskip 5pt

Assume that $G$ is  $2$-arc-transitive on some connected graph $\Ga=(V,E)$, and $\{\a,\b\}\in E$. Put $G^*=\l G_\a,G_\b\r$, the subgroup of $G$ generated by $G_\a\cup G_\b$. Then $|G:G^*|\leqslant 2$ with the equality holds if and only if $\Ga$ is bipartite and $G^*$ is the bipartition preserving subgroup of $G$, refer to \cite{Weiss}. Assume further that $\Ga$ is not a complete bipartite graph, and every minimal normal subgroup of $G$ contained in $G^*$ acts transitively on each of $G^*$-orbits on $V$. In 1993,
Praeger \cite{Praeger-ON, Praeger-ON-bi} proved that, except for one case when $\Ga$ is a bipartite graph, $G^*$ is a quasiprimitive permutation group of type HA, TW, AS or PA on each of its orbits, refer to \cite[Theorem 2]{Praeger-ON}, \cite[Theorems 2.1 and 2.3]{Praeger-ON-bi} and \cite[Theorem 6.1]{Prag-o'Nan-survey}. (Recall that a permutation group $G$ is quasiprimitive if every minimal normal subgroup of $G$ is transitive.) Roughly stated, either $(G,\Ga)$ is described as in \cite[Theorem 2.1 (c)]{Praeger-ON-bi}, or $G^*$ has a unique minimal normal subgroup  say $M$,
%the socle $\soc(G^*)$ of $G^*$,
and one of the following four cases occurs for  $M$ (and $G^*$):
\begin{itemize}%\itemindent=5pt
  \item[{\rm HA}] ({\em Holomorph Affine}): $M$ is abelian;

  \item[{\rm TW}] ({\em Twisted Wreath product}): $M$ is nonabelian and regular on each of $G^*$-orbits;

  \item[{\rm AS}] ({\em Almost Simple group}): $M$ is a nonabelian simple group;

  \item[{\rm PA}] ({\em Product Action}): $M=T_1\times\cdots\times T_n$ for some integer $n\geqslant 2$ and isomorphic nonabelian simple groups $T_i$, and for $\a\in V$ there are isomorphic subgroups $1\ne R_i<T_i$ such that $M_\a$ is a subdirect product of $R_1\times\cdots\times R_n$, that is, $M_\a$ projects surjectively onto every $R_i$.
\end{itemize}

For convenience, we say a connected $(G,2)$-arc-transitive  graph $\Ga$ is of HA, TW, AS or PA type if the case  HA, TW, AS or PA holds for $M$ and $G^*$, respectively. In addition, according to \cite{LS-pa},
 the type PA is said to be diagonal if each of the projections $M_\a\rightarrow R_i$ is injective, and non-diagonal otherwise.

\vskip 5pt

After Praeger's work, the existence of $2$-arc-transitive non-bipartite graphs with  HA, TW or AS  type was confirmed in just a few years. For example, the classification for graphs with HA type was given in \cite{Ivan-Pra}, constructions and examples  of graphs with TW  type were given in \cite{Baddeley} and \cite[Section 6]{Praeger-ON}, and some classification results of graphs with AS  type were given in \cite{FP1,FP2,HNP}. The existence problem of graphs with PA type was not answered until 2006 when Li and Seress \cite{LS-pa} constructed five families of $2$-arc-transitive non-bipartite graphs, four of them consist of graphs with diagonal PA type, and the other one consists of graphs of valency $9$ with non-diagonal PA type.

\vskip 5pt

In this paper, we first discuss some further structural information about $2$-arc-transitive (non-bipartite and bipartite) graphs with PA type. The following result is proved in Section \ref{proof}, which is helpful for us  to understand the behavior of $M_\a$ in the product action of  a $2$-arc-transitive group on some connected graph.

\begin{theorem}\label{th-1}
Let $\Ga=(V,E)$ be a connected $(G,2)$-arc-transitive graph with {\rm PA} type, and let $M=T_1\times\cdots\times T_n$, $G^*$ and $R_1$  be defined as above. Then, for $\a\in V$, one of the following holds.
\begin{itemize}
\item[(1)] %$M_\a\cong R_i$ for $1\leqslant i\leqslant n$.
$\Ga$ is of diagonal PA type.
   \item[(2)]  $M_\a\cong (\ZZ_p^k\times \ZZ_{m_1}){.}\ZZ_m$, $|\Ga(\a)|=p^k$, and $|R_1|$ is indivisible by $p^k$, where $m_1\div m$, $m\div (p^d-1)$ for some divisor $d$ of $k$ with $d<k$; in addition,
    \begin{itemize}

  \item[(i)] $n$ is divisible by some prime $r$, where either $r$ is an arbitrary  primitive prime divisor of $p^k-1$, or $(p,k)=(2,6)$ and $r\in \{3,7\}$; or
 \item[(ii)] $(p,k)=(2,6)$, and $M$ acts regularly on either the edge set or the arc set of $\Ga$; or
  \item[(iii)] $k=2$, and $p$ is a Mersenne prime.
  \end{itemize}
  \end{itemize}
\end{theorem}

\vskip 5pt

Li and Seress \cite{LS-pa} proved that, employing an equidistant linear $[4,2]_3$ code (see Section \ref{code} for the definition), one can construct $2$-arc-transitive graphs of valency $9$ with non-diagonal PA type from connected cubic graphs which admit a simple $2$-arc-regular group. This motivates us to develop a broader  construction for graphs with non-diagonal PA type. In Section \ref{code}, we
confirm that, for some suitable prime power $q$, there exist
equidistant  linear $[q+1,2]_q$ codes which admit a cyclic group of order $q^2-1$ acting regularly on the set of nonzero codewords. This allows us to construct some qusiprimitive permutation groups of PA type with a point stabilizer isomorphic to the affine group $\AGL_1(q^2)$, and then give a construction for $2$-arc-transitive graphs  with non-diagonal PA type.
Thus, in  Section \ref{sect=exam},  we construct  some $2$-arc-transitive graphs of valency $q^2$ with non-diagonal PA type, which meet Theorem \ref{th-1}~(i)  or (iii). Then, combining \cite[Lemma 5.2 and Example 5.3]{LS-pa}, we have the following result.

\begin{theorem}\label{th-2}
Let $q\geqslant 3$ be a prime power. Assume that $q+1$ has at most two distinct prime divisors, and either $q$ is even or $q\equiv -1\,(\mod 4)$.
Then there are  connected $2$-arc-transitive graphs of valency $q^2$ with non-diagonal {\rm PA} type.
\end{theorem}

We also construct  in  Section \ref{sect=exam} some graphs of valency $2^6$ and order
$2^{57}\cdot 3^{42}\cdot 7^{21}$, which give   examples for Theorem \ref{th-1}~(ii), see Example \ref{2^6}.

\vskip 5pt

For a graph $\Sigma=(V_0,E_0)$, the standard double cover $\Sigma^{(2)}$ is defined as the bipartite graph with vertex set  $V_0\times\ZZ_2$ such that $(\a_0,0)$ and $(\b_0,1)$ are adjacent
if and only if $\{\a_0,\b_0\}\in E_0$.
It is well-known that $\Sigma^{(2)}$ is connected if and only if
$\Sigma$ is connected and non-bipartite.
Define \[\iota: V_0\times\ZZ_2 \rightarrow V_0\times\ZZ_2,\,
(\a_0,i)\mapsto (\a_0,i+1).\]
Then $\iota\in \Aut(\Sigma^{(2)})$.
We view   $\Aut(\Sigma)$
as a subgroup of $\Aut(\Sigma^{(2)})$ in the following way
\[(\a_0,i)^g=(\a_0^g,i),\, \a_0\in V_0,\, i\in \ZZ_2,\,g\in \Aut(\Sigma).\] Then $\Aut(\Sigma^{(2)})$ has a subgroup
$\Aut(\Sigma)\times\l \iota\r$.
Thus, if $\Sigma$ is   $(G_0,2)$-arc-transitive   (and of some type) then $\Sigma^{(2)}$ is a $(G_0\times\l \iota\r,2)$-arc-transitive graph (of the same type).

\vskip 5pt

%Note, if $\Ga$ is non-bipartite then  the standard
%double cover $\Ga^{(2)}$ (see Section
%\ref{biparte-exam} for the definition)
%of $\Ga$ is of the same type as the graph $\Ga$.
Employing  standard double covers of graphs, one can easily get some firsthand examples of  bipartite $2$-arc-transitive graphs with HA, TW, AS or PA type.
In Section \ref{biparte-exam}, we give a construction for $2$-arc-transitive bipartite graphs of  diagonal PA type, which are not  standard double covers. In particular, the following result holds.
\begin{theorem}\label{th-3}
Let $p\geqslant 5$ be a prime.
Then there are  connected $2$-arc-transitive bipartite graphs of valency $p$ with  diagonal {\rm PA} type, which are not standard double covers of any graph.
\end{theorem}

\vskip 20pt

\section{On locally arc-transitive graphs}

In this section and the next section, we make some preparation for the proof of Theorem \ref{th-1}.

\vskip 5pt

Let $\Ga=(V,E)$ be a  graph, and $G\leqslant \Aut(\Ga)$.
The graph $\Ga$ is said to be $G$-locally  arc-transitive or $G$-locally  primitive if for every $\a\in V$, the stabilizer $G_\a$   acts transitively or primitively  on  $\Ga(\a)$, respectively.

\vskip 10pt

Let $\Ga=(V,E)$ be a connected graph, $\{\a,\b\}\in E$, $G\leqslant \Aut(\Ga)$ and $G^*=\l G_\a,G_\b\r$.
Assume that $G_\a$ and $G_\b$ act transitively on $\Ga(\a)$ and $\Ga(\b)$, respectively.
Then
$G^*$ acts transitively on $E$, and $G^*$  acts transitively on $V$ if $\Ga$ is not bipartite, refer to \cite[Exercise 3.8]{Weiss}. If $\Ga$ is not bipartite then $|G^*:G_\a|=|V|=|G:G_\a|$, yielding $G=G^*$.
%; of course, $G_\gamma\leqslant G^*$
%for every $\gamma\in V$.
Suppose that $\Ga$ is bipartite with two parts, say $U$ and $W$. Then $G^*$ fixes and acts transitively on both $U$ and $W$. Without loss of generality, let $\a\in U$ and $|U|\geqslant |W|$. We have
\[2|G^*:G_\a|=2|U|\ge |V|\geqslant |G:G_\a|.\]
It follows that either $G=G^*$, or $|G:G^*|=2$ and $G$ is transitive on $V$. In particular, $G^*$ is the  bipartition preserving subgroup of $G$, and thus   $G_\gamma\leqslant G^*$ for every $\gamma\in V$.
Now let $\gamma\in V$ and, without loss of generality, we set
$\gamma=\a^x$ for some $x\in G^*$.
Then $\Ga(\gamma)=\Ga(\a)^x$ and $G_\gamma=G_\a^x$. This implies that $G_\gamma$ acts transitively on $\Ga(\gamma)$, and the action is primitive if and only if $G_\a$ acts primitively on $\Ga(\a)$.
In summary, we have the following lemma.

\begin{lemma}\label{tech-1}
Let $\Ga=(V,E)$ be a connected graph, $\{\a,\b\}\in E$, $G\leqslant \Aut(\Ga)$ and $G^*=\l G_\a,G_\b\r$. Assume that $G_\a$ and $G_\b$ act transitively on $\Ga(\a)$ and $\Ga(\b)$, respectively.
Then  $\Ga$ is $G^*$-locally arc-transitive, and $\Ga$ is $G^*$-locally primitive if and only if   $G_\a$ and $G_\b$ act primitively on $\Ga(\a)$ and $\Ga(\b)$, respectively. Moreover, either
\begin{itemize}
  \item[(1)] $\Ga$ is not bipartite, and $G=G^*$ is transitive on $V$; or
  \item[(2)]  $\Ga$ is a bipartite graph with two parts the $G^*$-orbits on $V$, and $|G:G^*|\leqslant 2$, where the equality holds if and only if $G$ is transitive on $V$.
\end{itemize}
\end{lemma}

For    locally primitive graphs, by \cite[Lemmas 2.5 and 2.6]{567},
the next result holds.

\begin{lemma}\label{tech-2}
Assume $\Ga=(V,E)$ is a connected $G$-locally primitive graph, and   $N$ is a normal subgroup of $G$.
\begin{itemize}
  \item[(1)] If $G$ is transitive on $V$ and $N_\a\ne 1$ for some $\a\in V$ then $\Ga$ is $N$-locally arc-transitive.
  \item[(2)] If $N$ is intransitive on each of $G$-orbits on $V$, then either
      %$N$ is semiregular on $V$,
      %that is, $N_\b=1$ for all $\b\in V$.
  \begin{itemize}
    \item[(i)] $N$ is semiregular on $V$,
      that is, $N_\a=1$ for all $\a\in V$, and $N$ itself is the kernel of $G^*$ acting on the $N$-orbits; or
    \item [(ii)] $G$ is transitive on $V$, $N$ has two orbits on $V$, and either $N$ is semiregular on $V$ or $\Ga$ is $N$-locally arc-transitive.
  \end{itemize}
\end{itemize}
\end{lemma}

The next lemma says that some conclusion in Lemma \ref{tech-2}
is   true for a bipartite graph $\Ga$ under some
weaker conditions. For
$U_1,W_1\subseteq V$,  denote  by $[U_1,W_1]$ the subgraph of $\Ga$ induced by $U_1\cup W_1$.

\begin{lemma}\label{tech-3}
Let $\Ga=(V,E)$ be a connected bipartite graph, $\{\a,\b\}\in E$, $G\leqslant \Aut(\Ga)$ and $G^*=\l G_\a,G_\b\r$. Assume that $G_\a$ acts primitively on $\Ga(\a)$, and $G^*$ has a normal subgroup $N$ which is intransitive on each of $G^*$-orbits on $V$. Then   $N$ is semiregular on $V$, and $N$ itself is the kernel of $G^*$ acting on the $N$-orbits.
\end{lemma}
\proof
Let $U$ and $W$ be the $G^*$-orbits containing $\a$ and $\b$, respectively. For an arbitrary $\gamma\in U$,  we have
$\gamma=\a^x$ for some $x\in G^*$, and thus
 $\Ga(\gamma)=\Ga(\a)^x$ and $G_\gamma=G_\a^x$, it follows that $G_\gamma$ acts primitively on $\Ga(\gamma)$.

 \vskip 5pt

Let $\mathcal{U}$ and $\mathcal{W}$ be the sets of $N$-orbits on $U$ and $W$, respectively. Pick $U_1\in \mathcal{U}$  and $\gamma\in U_1$.  Then $\{\Ga(\gamma)\cap W_1\mid  W_1\in \mathcal{W},\, \Ga(\gamma)\cap W_1\ne\emptyset\}$ is a $G_\gamma$-invariant partition of $\Ga(\gamma)$. Since $G_\gamma$ acts primitively on $\Ga(\gamma)$, either $\Ga(\gamma)\subseteq  W_1$ for some $W_1\in \mathcal{W}$, or  $[U_1,W_1]$ is a matching without isolated vertex for every $W_1\in \mathcal{W}$ with $\Ga(\gamma)\cap W_1\ne\emptyset$.

  \vskip 5pt

 Suppose first that $\Ga(\gamma)\subseteq  W_1$ for some $W_1\in \mathcal{W}$. Then every vertex in $U_1$ has no neighbor in $W\setminus W_1$ and, since  $W_1$ is an $N$-orbit, every vertex in $W_1$ has some neighbor in $U_1$.
Let $\delta\in W_1$, and pick its neighbors $\gamma_1$ and $\gamma_2$ with $\gamma_1\in U_1$. Let $U_2$ be the $N$-orbit containing $\gamma_2$. Then $U_1^y=U_2$, where $y\in G_\delta$ with $\gamma_1^y=\gamma_2$. Noting that $W_1^y=W_1$, it follows that  $[U_1,W_1]$   and    $[U_2,W_1]$    are isomorphic. Thus every vertex in $U_2$ has no neighbor in $W\setminus W_1$. Let $U_0$ be the set of vertices which have neighbors in $W_1$. By the above argument,  every vertex in $U_0$ has no neighbor in $W\setminus W_1$ and, by the choice of $U_0$, every vertex in $W_1$ has no neighbor in $U\setminus U_0$.
It follows that $\Ga=[U_0,W_1]$, and then $W_1=W$, which contradicts that $N$ is intransitive on $W$.

 \vskip 5pt

Now, for arbitrary $U_1\in \mathcal{U}$ and $W_1\in \mathcal{W}$, the subgraph $[U_1,W_1]$ is either a empty graph or a matching without isolated vertex. Let $K$ be the kernel of $G^*$ acting on $\mathcal{U}\cup \mathcal{W}$. We have $N\leqslant K$.
In the following, we will show that $K_\gamma=1$ for all $\gamma\in V$,
and then the lemma follows.

 \vskip 5pt

Let $\gamma,\,\delta\in V$. Since $\Ga$ is connected, pick a path $\gamma=\a_0,\a_1,\ldots, \a_n=\delta$ from $\gamma$ to $\delta$.
For $0\leqslant i\leqslant n$, let $V_i$ be the $N$-orbit containing $\a_i$. Suppose that $K_\gamma$ fixes $\a_{i-1}$. Noting that $K_\gamma$ fixes both $V_{i-1}$ and $V_i$ set-wise, since $\a_{i-1}$ has a unique neighbor in $V_i$, it follows that
$K_\gamma\leqslant K_{\a_i}$. By induction, we have $K_\gamma\leqslant K_{\delta}$. Thus $K_\gamma$ fixes $V$ point-wise, and hence $K_\gamma=1$. This completes the proof.
\qed

 \vskip 5pt

\begin{lemma}\label{tech-5}
Let $\Ga=(V,E)$ be a connected $G$-locally arc-transitive graph, $\{\a,\b\}\in E$ and $N\unlhd G$. Suppose that $(|N_\a|,|\Ga(\a)|)=1=(|N_\b|,|\Ga(\b)|)$. Then $N$ is semiregular on $V$.
\end{lemma}
\proof
Let $\gamma$ be an arbitrary vertex of $\Ga$. By the assumption, since $G$ acts transitively on $E$, we have  $(|N_\gamma|,|\Ga(\gamma)|)=1$.
Note that $N_\gamma\unlhd G_\gamma$ and $G_\gamma$ acts transitively on $\Ga(\gamma)$. Then all $N_\gamma$-orbits on $\Ga(\gamma)$ have the same length, which is a common divisor of $|\Ga(\gamma)|$ and $|N_\gamma|$. It follows that $N_\gamma$ fixes $\Ga(\gamma)$ point-wise. In particular, $N_\gamma\leqslant N_\delta$ for $\delta\in \Ga(\gamma)$. Again since $(|N_\delta|,|\Ga(\delta)|)=1$, a similar argument implies that
$N_\delta$ fixes $\Ga(\delta)$ point-wise, and so $N_\gamma$ fixes $\Ga(\delta)$ point-wise. Thus, since $\Ga$ is connected, we conclude that $N_\gamma$ fixes $V$ point-wise, and so $N_\gamma=1$. Then $N$ is semiregular on $V$.
\qed

\vskip 20pt

\section{Two elementary results on primitive affine groups}
Recall that, for positive integers $p, k>1$, a primitive prime divisor of $p^k-1$ is a prime which divides $p^k-1$  but does not divide $p^i-1$ for all $0< i<k$. If $r$ is a  primitive prime divisor of $p^k-1$, then $k$ is the smallest positive integer
with $p^k\equiv 1\,(\mod r)$,
 and thus $k$ is a divisor of $r-1$; if further $r\div (q^l-1)$ with $l\geqslant 1$ then $k\div l$.
These facts yield a criterion for  affine  primitive permutation groups.

 \vskip 5pt
For a group $X$ and subgroups $Y,Z\leqslant X$, let $\bC_Y(Z)=\{y\in Y\mid yz=zy \mbox{ for all }z\in Z\}$, called the centralizer of $Z$ in $Y$.

\begin{lemma}\label{tech-4}
Let $H$ be a  permutation group  on a set $\Omega$, and $\a\in \Omega$. Suppose that $H$ has a regular  normal subgroup $P\cong \ZZ_p^k$,   where $k\geqslant 2$ and $p$ is a prime. Suppose that $p^k-1$ has a  primitive prime divisor $r$, and $|H_\a|$ is divisible by $r$. Then $H$ is primitive on $\Omega$.
\end{lemma}
\proof
Let $Q$ be a Sylow $r$-subgroup of $H_\a$. Then $Q\ne 1$ as $r$ is a divisor of $|H_\a|$. Set $K=PQ$. We next show that $K$ is primitive on $\Omega$. It suffices to prove that  $Q$ is a maximal subgroup of $K$.

 \vskip 5pt

%It suffices to show that $Y$ is primitive on $\Omega$, i.e,
%$R$ is a maximal subgroup of $Y$.
By Maschke's Theorem (refer to \cite[p.123, I.17.7]{Huppert}), since $(p,|Q|)=1$,
we have $P=P_1\times\cdots\times P_l$, where $P_i$ are minimal $Q$-invariant subgroups of $P$. Considering the conjugation of $Q$ on $P_i$, the group $Q$ induces a subgroup of the automorphism group $\Aut(P_i)$ of $P_i$ with kernel  $\bC_Q(P_i)$. Let $|N_i|=p^{k_i}$. Then
$\Aut(P_i)$ is isomorphic to the general linear group $ \GL_{k_i}(p)$, and so
\[Q/\bC_Q(P_i)\lesssim \Aut(P_i)\cong \GL_{k_i}(p),\, 1\leqslant i\leqslant l.\]

\vskip 5pt

Suppose that $l>1$. Then $k_i<k$ for every $i$, and so $|\GL_{k_i}(p)|$ is indivisible by $r$. It follows that $Q=\bC_Q(P_i)$ for all $i$, and thus $Q$ centralizes $P$. Then $Q\unlhd K$, which is impossible as $1\ne Q=K_\a$.
%. Since $K$ is transitive on $\Omega$, every $R$-orbit has
%a length a divisor of $p^k$. On the other hand, the length
%of every $R$-orbit is a power of $r$. Then $R$ fixes $\Omega$
%point-wise, and so $R=1$, a contradiction.
Therefore, $l=1$, which yields that $P$ is a minimal normal subgroup of $K$.

 \vskip 5pt

Let $L$ be a maximal subgroup of $K$ with $Q\leqslant L$. Then
$K>L=PQ\cap L=(P\cap L)Q$, and so $P\cap L\ne P$. Since $P$ is abelian and $P\unlhd K$, we have $P\cap L\unlhd P$ and $P\cap L\unlhd L$, and thus $P\cap L\unlhd \l P,L\r=K$. Then $P\cap L=1$ as  $P$ is a minimal normal subgroup of $K$. Thus $L=(P\cap L)Q=Q$. This says that $Q$ is a maximal subgroup of $K$, and then $K$ is primitive on $\Omega$. Noting that $K\leqslant H$, the lemma follows.
\qed

 \vskip 5pt

A transitive permutation group $H$  on a set $\Omega$ is  a Frobenius group
if $H_\a\ne 1$  for $\a\in \Omega$, and $H_{\a\b}=1$ for all $\beta\in \Omega\setminus \{\a\}$. The following lemma gives a characterization of imprimitive Frobenius groups with abelian socle, see \cite[Lemma 2.2]{basic-two} for example.
Recall that, for a finite group $X$, the socle $\soc(X)$ of $X$ is generated by all minimal normal subgroups of $X$.

 \vskip 5pt

\begin{lemma}\label{Frobenius}
Let $K$ be an imprimitive Frobenius group on $\Omega$ with $\soc(K)=P\cong \ZZ_p^k$, where $p$ is a prime  and   $k\geqslant 2$. Then $K_\a$ is isomorphic to an irreducible subgroup of the general linear group $\GL_l(p)$ for some $l$, and $|K_\a|$ is a divisor of $p^d-1$, where $2l\leqslant k$ and
$d$ is a common divisor of $k$ and $l$.
 \end{lemma}

\vskip 5pt

\begin{lemma}\label{tech-6}
Let $H$ be a  $2$-transitive affine group  of degree $2^6$ on a set $\Omega$, and let $1\ne K\unlhd H$.
Assume that $K_\a\ne 1$ for $\a\in V$, and $K$ is imprimitive on $\Omega$. Then
\begin{itemize}
\item[(1)] $K_\a\cong \ZZ_s$ with $s\in \{3,7\}$, and there is $x\in H_\a$   such that $K_\a\l x\r\cong \ZZ_{21}$; and
\item[(2)] for each $x\in H_\a$ with $K_\a\l x\r\cong \ZZ_{21}$, the subgroup $K\l x\r$ is  primitive  on $\Omega$.
\end{itemize}
\end{lemma}
\proof
By  \cite[pp.215-217, Theorems 7.2C and 7.2E]{Dixon}, $K$ is an imprimitive Frobenius group. Applying Lemma \ref{Frobenius}, we get $K_\a\cong \ZZ_3$ or $\ZZ_7$. Calculation with GAP \cite{GAP} shows that there are eleven  $2$-transitive affine groups  of degree $2^6$ containing an imprimitive Frobenius subgroup. Checking one by one these groups, we conclude that
either $K_\a\cong \ZZ_3$ is the center of $H_\a$, or $K_\a$ is contained in a cyclic subgroup of order $21$ in $H_\a$. Then part (1) of this lemma follows.

\vskip 5pt

Assume that $x\in H_\a$ with $K_\a\l x\r\cong \ZZ_{21}$, and set $X=K\l x\r$. Then $\soc(H)\unlhd X$ and $X_\a\cong \ZZ_{21}$. Without loss of generality, we assume that $K_\a\cap \l x\r=1$, let $ \l x\r\cong \ZZ_r$ and write $X_\a=\l y\r\times \l x\r$ with $K_\a=\l y\r\cong \ZZ_s$.

\vskip 5pt

By Maschke's Theorem, we have
$\ZZ_2^6\cong \soc(H)=P_1\times\cdots\times P_l$, where $P_i$ are minimal $X_\a$-invariant subgroup of $\soc(H)$. Since $K$ is an imprimitive Frobenius group, $y$ does not centralizes every $P_i$, and  $s$ is a divisor of $|P_i|-1$, refer to \cite[p.191, (35.25)]{Aschbacher}.
Suppose that $l>1$. Then either
$s=3$, $P_i\cong \ZZ_2^2$ and  $l=3$, or $s=7$, $P_i\cong \ZZ_2^3$ and  $l=2$, where $1\leqslant i\leqslant l$. Note that $X_\a/\bC_{X_\a}(P_i)\lesssim \Aut(P_i)$. Assume first that $s=3$. Then $r=7$, and $\Aut(P_i)\cong \GL_2(2)\cong\Sy_3$. This implies that $x$ centralizes every $P_i$. Thus $\l x\r\unlhd H$, which is impossible as $1\ne \l x\r\leqslant X_\a$.
Now let $s=7$. Then $l=2$, and $P_1\cong P_2\cong \ZZ_2^3$. We have $\l y\r\cong (\l y\r\bC_{X_\a}(P_i))/\bC_{X_\a}(P_i) \leqslant X_\a/\bC_{X_\a}(P_i)\lesssim\Aut(P_i)\cong \GL_3(2)$. By the Atlas \cite{Atlas}, $\GL_3(2)$ has no element of order $21$. It follows that $x$ centralizes every $P_i$, which leads to a similar contradiction as above. Therefore, $l=1$, and then $\soc(H)$ is a minimal normal subgroup of $X$.  Thus $X_\a$ is a maximal subgroup of $X$, and   part (2) of this lemma follows.
\qed

\vskip 20pt

\section{The proof of Theorem \ref{th-1}}\label{proof}

Let $\Ga=(V,E)$ be connected graph of valency no less than $3$, and $G\leqslant \Aut(\Ga)$. Let   $G^*=\l G_{\a_1},G_{\a_2}\r$ for some $\{\a_1,\a_2\}\in E$, and let $M=\soc(G^*)$. Assume that $\Ga$ is $(G,2)$-arc-transitive, and   $G^*$ is a quasiprimitive group of PA type on each of   $G^*$-orbits.  Then both $G^*$ and $M$ have the same orbits on $V$.
By \cite{Praeger-ON,Praeger-ON-bi}, we have
\begin{itemize}
  \item[(I)] $M=T_1\times T_2\times \cdots \times T_n$ is the unique minimal normal subgroup of $G^*$, where $n\geqslant 2$ and $T_i$ are isomorphic nonabelian simple groups; and
  \item[(II)] for $\a\in V$, there are   subgroups $R_i<T_i$ such that  $M_\a\leqslant R_1\times \cdots \times R_n$ and, for every $i$, the projection
\[\pi_i: M_\a\rightarrow R_i,\,\, x_1x_2\cdots x_n\mapsto x_i,\, \mbox{ where }x_j\in R_j \mbox{ for all } j\]
is a surjective group homomorphism.
\end{itemize}

 \vskip 5pt

Note that  $T_1,\,T_2,\,\ldots,\,T_n$ are all minimal normal subgroups of $M$, refer to \cite[p.51, I.9.12]{Huppert}. Since $M$ is a minimal normal subgroup of $G^*$, we have  
\begin{itemize}
  \item[(III)]
$G_\a$ acts transitively on $\{T_1,T_2,\ldots,T_n\}$ by conjugation.
\end{itemize}
Clearly, $M_\a\unlhd G_\a$. For $h\in G_\a$, letting $T_i^h=T_{i'}$, we have
\[R_{i'}=\pi_{i'}(M_\a)=\pi_{i'}(M_\a^h)\leqslant \pi_i(M_\a)^h=R_i^h.\]
It follows that
\begin{itemize}
  \item[(IV)]
$G_\a$ acts transitively on $\{R_1,R_2,\ldots,R_n\}$ by conjugation; in particular, $R_1\cong \cdots\cong R_n$.
\end{itemize}

\vskip 5pt

\vskip 10pt

 For convenience, we set $N_i=\prod_{j\ne i}T_j$, where $1\leqslant i\leqslant n$. Then
 \begin{itemize}
  \item[(V)]
 $N_i\lhd M$, and the kernel $\ker(\pi_i)$ of $\pi_i$ equals to $(N_i)_\a$.
 \end{itemize}
 Note that $N_1$, $\ldots,\,N_n$ are all maximal normal subgroups of $M$, refer to \cite[p.51, I.9.12]{Huppert}. We have
 \begin{itemize}
  \item[(VI)]
$G_\a$ acts transitively on both $\{N_1,N_2,\ldots,N_n\}$ and $\{\ker(\pi_1),\ldots, \ker(\pi_n)\}$ by conjugation; in particular, $\ker(\pi_1)\cong \cdots\cong \ker(\pi_n)$.
\end{itemize}
In addition,   the following lemma holds.

\begin{lemma}\label{ker-pi}
 Every $N_i$ is intransitive on each of   $M$-orbits  on $V$.
\end{lemma}
\proof  Suppose that some
 $N_i$ acts transitively on one of the $M$-orbits.
Then $M=N_iM_\gamma$ for some $\gamma\in V$. Thus
$T_i\cong M/N_i=N_iM_\gamma/N_i\cong  M_\gamma/(N_i)_\gamma$.
Then $M_\gamma$ has a composition factor isomorphic to $T_i$, which is impossible as $M_\gamma\cong M_\a\leqslant R_1\times \cdots \times R_n$. This completes the proof.
\qed

\vskip 10pt

In the following, we will formulate the case where
some   $\pi_i$ is not injective.

\vskip 5pt

By Lemma \ref{tech-2},   $\Ga$ is $M$-locally arc-transitive.
If $\Ga$ is $M$-locally primitive, then Theorem \ref{th-1} is true by the following simple lemma.

\vskip 5pt

\begin{lemma}\label{local-prim} Assume $\Ga$ is $M$-locally primitive. Then every $\pi_i$ is injective; in particular, $M_\a\cong R_i$ for all $i$.
\end{lemma}
\proof
Suppose that some $\pi_i$  is not injective. Then $\pi_i$ has nontrivial kernel $\ker(\pi_i)=(N_i)_\a$.
Then, by Lemmas \ref{tech-1} and \ref{tech-2}, $N_i$ is transitive on one of the $M$-orbits on $V$, which contradicts Lemma \ref{ker-pi}. This completes the proof.
\qed

\vskip 10pt

We next deal with the case where $\Ga$ is not   $M$-locally primitive.  For $X\leqslant G$, denote by $X_\a^{[1]}$ the kernel of $X_\a$ acting on $\Ga(\a)$, and by $X_\a^{\Ga(\a)}$ the permutation group induced by  $X_\a$  on $\Ga(\a)$.
By \cite{basic-two}, we have the following lemma.

\begin{lemma}\label{imp-stab}
If  $\Ga$ is not $M$-locally primitive, then 
one of the following holds.
 \begin{itemize}
          \item[(1)] $M_\a\cong (\ZZ_p^k\times \ZZ_{m_1}){.}\ZZ_m$, $|\Ga(\a)|=p^k$ and $M_\a^{[1]}\cong \ZZ_{m_1}$, where $m_1\div m$, $m\div (p^d-1)$ for some divisor $d$ of $k$ with $d<k$;
        \item[(2)] $M_\a\cong (\ZZ_3^4\times Q){.}\Q_8$, $|\Ga(\a)|=3^4$ and $M_{\a}^{[1]}\cong Q$, where $Q$ is isomorphic to a subgroup of the quaternion group $\Q_8$.
\end{itemize}
\end{lemma}

\vskip 10pt

Together with Lemmas \ref{local-prim} and \ref{imp-stab}, the following lemma fulfills the proof of Theorem \ref{th-1}.

\begin{lemma}\label{prim-div}
Assume that $|\Ga(\a)|=p^k$ and $M_\a$ is described as in {\rm (1)} or {\rm (2)} of Lemma {\rm \ref{imp-stab}}. Let $p^l$ be the highest power of $p$ dividing $|R_1|$.
\begin{itemize}
\item[(1)] If $l=k$ then every $\pi_i$ is injective.
   \item[(2)]  If $l<k$ then one of the follows holds.
    \begin{itemize}

  \item[(i)] $n$ is divisible by some prime $r$, where either $r$ is an arbitrary  primitive prime divisor of $p^k-1$, or $(p,k)=(2,6)$ and $r\in \{3,7\}$;
 \item[(ii)] $(p,k)=(2,6)$, and $M$ acts regularly on the edge set or arc set of $\Ga$;
  \item[(iii)] $k=2$, and $p$ is a Mersenne prime.
  \end{itemize}
  \end{itemize}
\end{lemma}
\proof
Recalling that $\pi_1:M_\a\rightarrow R_1$ is a surjective homomorphism, we have $l\leqslant k$.
Assume that $l=k$. Then every $\ker(\pi_i)$ has order indivisible by $p$. Noting $(N_i)_\a=\ker(\pi_i)$, by Lemma \ref{tech-5}, $\ker(\pi_i)=1$, and  part (1) of this lemma is true.

 \vskip 5pt

Assume that  $l<k$ from now on.
If $p^k-1$ has no primitive prime divisor then, by Zsigmondy's Theorem,  either $(p,k)=(2,6)$, or $k=2$ and $p$ is a Mersenne prime. The latter case is just the case (iii) of the lemma. For
 $(p,k)=(2,6)$, if $M_\a\cong \ZZ_2^6$ then we get the case (ii) of this lemma.

 \vskip 5pt

In the following, we   assume  further that either $(p,k)=(2,6)$ and $M_\a\not\cong \ZZ_2^6$, or $p^k-1$ has a primitive prime divisor $r$.
Noting that $G_\a$ acts $2$-transitively on $\Ga(\a)$, it follows that  $p^k-1$ is a divisor of $|G_{\a\b}|$ for $\b\in \Ga(\a)$, and then either $21$ or $r$ is a divisor of $|G_{\a\b}|$, respectively.
In addition, for $(p,k)=(2,6)$, we have $M_\a^{\Gamma(\a)}\cong \ZZ_2^6{:}\ZZ_s$ with $s\in \{3,7\}$ by Lemma \ref{imp-stab}; in this case, we set  $r={21\over s}$.

 \vskip 5pt

 {\it Claim }1. If $(p,k)=(2,6)$ then there is an element $x\in G_{\a\b}$ of order $r$ such that $M_{\a\b}\l x\r=M_{\a\b}\times \l x\r$, where $\b \in \Ga(\a)$.

  \vskip 5pt

  Assume that $(p,k)=(2,6)$. By  Lemma \ref{imp-stab}, we conclude that $M_{\a\b}$ is an abelian group of order $s$ or $s^2$. Then $M_{\a\b}\cong \ZZ_s$, $\ZZ_s^2$ or $\ZZ_{s^2}$, and thus $\Aut(M_{\a\b})$ has order $s-1$, $s(s-1)(s^2-1)$ or $s(s-1)$, respectively. Since $M_{\a\b}\unlhd G_{\a\b}$, every element in $G_{\a\b}$ induces an automorphism of $M_{\a\b}$ by conjugation. If $s=3$ then  $|\Aut(M_{\a\b})|$ is indivisible by $r=7$, and so $M_{\a\b}$ is centralized by every element  of order $7$ in $G_{\a\b}$, our claim is true in this case.

   \vskip 5pt

  Now let $s=7$ and $r=3$. Then a Sylow $3$-subgroup of $\Aut(M_{\a\b})$ is isomorphic to $\ZZ_3$, $\ZZ_3^2$ or $\ZZ_3$ when $M_{\a\b}\cong \ZZ_7$, $\ZZ_7^2$ or $\ZZ_{7^2}$, respectively. Noting that the $2$-transitive affine group $G_\a^{\Ga(\a)}$ has a normal subgroup isomorphic to $\ZZ_2^6{:}\ZZ_7$, calculation with GAP \cite{GAP} shows that $(G_\a^{\Ga(\a)})_\b$ has a subgroup isomorphic to $\ZZ_9$. Pick a Sylow $3$-subgroup $Q$ of $G_{\a\b}$. Then $Q$ acts unfaithfully on $M_{\a\b}$ by conjugation; otherwise, $Q\lesssim\ZZ_3^2$, which is impossible. Thus $\bC_Q(M_{\a\b})\ne 1$, and every element of order $3$ in $\bC_Q(M_{\a\b})$  is a desired $x$. Then Claim 1 follows.

\vskip 5pt

Now fix  an element $x\in G_{\a\b}$ of order $r$, where either $r$ is a  primitive prime divisor of $p^k-1$, or $(p,k)=(2,6)$, $r={21\over s}$ and $x$ is described as in Claim 1. Then $M\cap \l x\r=M_\a\cap \l x\r=1$.
Set $X=M\l x\r$. Clearly, $\Ga$ is $X$-locally arc-transitive, and $|X_\gamma|=r|M_\a|$ for all  $\gamma\in V$. In addition, for $(p,k)=(2,6)$, we have $X_{\a\b}=M_{\a\b}\times \l x\r$.

 \vskip 5pt

{\it Claim} 2. Either $X_\a$ acts primitively on $\Ga(\a)$, or $X_\b$ acts primitively on $\Ga(\b)$.

 \vskip 5pt

By Lemma \ref{imp-stab}, either $M_\a\cong \ZZ_2^2$ or $|\Ga(\a)|\geqslant 8$. Assume first $M_\a\cong \ZZ_2^2$.
Then $r=3$, $X_\a=M_\a\l x\r$, $X_\b=M_\b\l x\r$ and $X_{\a\b}=\l x\r$. Suppose that $X_\a^{[1]}\ne 1 \ne X_\b^{[1]}$. Then $X_\a^{[1]}=X_\b^{[1]}=X_{\a\b}=\l x\r$, yielding $\l x\r\unlhd \l X_\a,X_\b\r$. Note that $\l X_\a,X_\b\r$ acts transitively on $E$, refer to \cite[Exercise 3.8]{Weiss}. It follows that $\l x\r$ fixes every edge of $\Ga$, and thus $\l x\r=1$, a contradiction. We have
$X_\a^{[1]}=1$ or $X_\b^{[1]}=1$. Then one of $X_\a^{\Ga(\a)}$ and $X_\b^{\Ga(\b)}$ is a $2$-transitive group of degree $4$, and Claim 2 is true in this case.

 \vskip 5pt

Assume that $|\Ga(\a)|\geqslant 8$. Then, by \cite[Theorem 4.7]{Weiss}, $G_\a^{[1]}\cap G_\b^{[1]}=1$, and so $X_\a^{[1]}\cap X_\b^{[1]}=1$. Considering the actions of $X_{\a\b}$ on $\Ga(\a)$ and $\Ga(\b)$, we have
\[X_{\a\b}^{\Ga(\a)}\cong X_{\a\b}/X_\a^{[1]}, \, X_{\a\b}^{\Ga(\b)}\cong X_{\a\b}/X_\b^{[1]}.\]
If neither $|X_{\a\b}^{\Ga(\a)}|$ nor $|X_{\a\b}^{\Ga(\b)}|$ is divisible by $r$, then all Sylow $r$-subgroups are contained in both $X_\a^{[1]}$ and  $X_\b^{[1]}$, which contradicts that $X_\a^{[1]}\cap X_\b^{[1]}=1$. Without loss of generality, we
assume that $|X_{\a\b}^{\Ga(\a)}|$ is divisible by $r$.
If $r$ is a primitive prime divisor of $p^k-1$, then $X_\a^{\Ga(\a)}$ is primitive by Lemma \ref{tech-4}.
Now let $(p,k)=(2,6)$. Noting that $\ZZ_s\cong M_{\a\b}^{\Ga(\a)}\unlhd X_{\a\b}^{\Ga(\a)}$ and $X_{\a\b}\l x\r=M_{\a\b}\times \l x\r$, we have $X_{\a\b}^{\Ga(\a)}\cong \ZZ_{21}$. Then $X_\a^{\Ga(\a)}$ is primitive by Lemma \ref{tech-6}. Thus Claim 2 follows.

\vskip 5pt

Finally, consider the action of $\l x\r$ on $\{T_1,\ldots,T_n\}$ by conjugation.
Suppose  that some $T_i$, say $T_1$ without loss of generality, is normalized by $x$. Then $N_1=\prod_{j\ne 1}T_j$ is also  normalized by $x$, and thus $N_1\unlhd X$. Note that $N_1$ is intransitive on each of $M$-orbits, see Lemma \ref{ker-pi}. Assume that $\Ga$ is not bipartite. Then, by  Claim 2 and Lemma \ref{tech-2}, $N_1$ is semiregular on $V$, and so $\ker(\pi_1)=(N_1)_\a=1$, yielding
$M_\a\cong R_1$. Thus $k=l$, which is not the case. If $\Ga$ is   bipartite then, by Claim 2 and Lemma \ref{tech-3}, $N_1$ is semiregular on $V$, we have a similar contradiction as above.
Therefore, $\l x\r$ acts faithfully and semiregularly on $\{T_1,\ldots,T_n\}$.
Then $r$ is a divisor of $n$, and case (i) of this lemma follows.
This completes the proof.
\qed

\vskip 20pt

\section{A construction of equidistant linear codes}\label{code}
Let $q=p^f$ for some prime $p$ and integer $f\geqslant 1$. Denote by $\FF_q$ the field of order $q$, and $\FF_q^n$ the $n$-dimensional row vector space over $\FF_q$, where $n\geqslant 1$. For a vector $\bv=(v_1,v_2,\ldots, v_n)\in \FF_q^n$, letting $\supp(\bv)=\{i \mid v_i\ne 0,\,1\leqslant i\leqslant n\}$,
the weight $\wt(\bv)$ is defined as $|\supp(\bv)|$, i.e.,
the number of nonzero coordinates of $\bv$.

\vskip 5pt

Let $k$ be an integer with $1\leqslant k\leqslant n$. Every $k$-dimensional subspace $\cC$ of $\FF_q^n$ is called a linear  $[n,k]_q$ code, where $n$ is called the length  of $\cC$, and the vectors in $\cC$ are called codewords.
A linear $[n,k]_q$ code  $\cC$ is said to be equidistant if all nonzero codewords have the same weight say $\omega$, while $\omega$ is called the weight of  $\cC$ and write $\wt(\cC)=\omega$.

\vskip 5pt

Let $\cC$ be an equidistant  linear $[n,2]_q$ code with $\wt(\cC)=\omega$.
For ${\bf 0}\ne \bw\in \cC$, define
\[\cC_{\bw}=\{\bu\in\cC\mid \supp(\bu)=\supp(\bw) \mbox{ or } \emptyset\}.\]
Then it is easily shown that $\cC_{\bw}$ is a $1$-dimensional subspace of $\cC$, and every $1$-dimensional subspace of $\cC$ is obtained in the form of $\cC_{\bw}$. Choose $\bw_\ell\in \cC$, $1\leqslant \ell\leqslant q+1$, with $\cC=\cup_{\ell=1}^{q+1}\cC_{\bw_\ell}$. Set $\Delta=\cup_{\ell=1}^{q+1}\supp(\bw_\ell)$, and view $\cC$ as an $[m,2]_q$ code, where $m=|\Delta|$. Then $\omega\leqslant m-1$ by the Singleton bound, refer to \cite[p.73, Corollary 2.50]{Hir}. Consider the linear maps $\pi_i: \cC\rightarrow \FF_q$ given by $(v_1,\ldots,v_n)\mapsto v_i$, where $i\in \Delta$. Clearly,  every $\pi_i$ is surjective, and
$\ker(\pi_i)$ is $1$-dimensional. Then, for each $i\in \Delta$, there is some $\bw_\ell$ with $i\not\in \supp(\bw_\ell)$ and $\ker(\pi_i)= \supp(\bw_\ell)$. It follows that $\Delta\setminus \supp(\bw_\ell)$, $1\leqslant \ell\leqslant q+1$, are disjoint subsets of $\Delta$. Noting that $|\Delta\setminus \supp(\bw_\ell)|=m-\omega$, we have
\[q+1\leqslant {m\over m-\omega}\leqslant {n\over n-\omega}.\]
 Then we get the following fact.

\begin{lemma}\label{eqid}
Let $\cC$ be an equidistant  linear $[n,2]_q$ code with $\wt(\cC)=\omega$. If $n=q+1$ then $\omega=n-1$, and $\ker(\pi_i)$, $1\leqslant i\leqslant n$, are distinct $1$-dimensional subspaces of $\cC$.
\end{lemma}

\vskip 5pt
%\vskip 20pt

Let $n=q +1$ from now on. Denote by  $\FF_q^*$   the multiplicative group of $\FF_q$, and write  $\FF_q^*=\l\eta,\lambda\r$, where $\lambda$ has odd order, and $\eta$ has order a power of $2$. Clearly, $\FF^*=\l\eta\lambda\r=\l\eta\lambda^2\r$. Note that $\eta=1$ if $q$ is even, and $(\eta\lambda)^{q-1\over 2}=-1=(\eta\lambda^2)^{q-1\over 2}$ if $q$ is odd. Pick two invertible $n\times n$ matrices over $\FF_q$:
\[ \bD=\left(
     \begin{array}{ccc}
      \eta\lambda & 0 & {\bf 0} \\
    0 & \lambda &  {\bf 0} \\
     0 & 0  &  \eta\lambda \bI_{n-2}
     \end{array}
   \right),\,\,\, \bP=\left(
                        \begin{array}{cc}
                          {\bf 0}' & \bI_{n-1} \\
                          1 & {\bf 0} \\
                        \end{array}
                      \right),
\]
where $\bI_m$ denotes the identity matrix of order $m$. Let $\bA=\bD\bP$. Then
\[\bA^n=\eta^{n-1}\lambda^n \bI_n=\eta\lambda^2\bI_n.\]
In particular, $\bA$ has order $n(q-1)$ as an element of the general linear group $\GL_n(q)$.

\vskip 5pt

View $\bA$ as the  linear transformation of $\FF_q^n$ given by right multiplication on the row vectors. Then we have an action of the cyclic group $\l \bA\r$ on $\FF_q^n$.
A linear $[n,k]_q$ code $\cC$ is said to be $\l\bA\r$-invariant if $\bu\bA\in \cC$ for all $\bu\in \cC$, and $\l\bA\r$-irreducible if further $\cC$ does not contains a $\l\bA\r$-invariant  linear $[n,k']_q$ code for some $1\leqslant k'<k$. A $\l\bA\r$-invariant linear  $[n,k]_q$ code $\cC$ is said to be faithful if $\l \bA\r$ acts faithfully on $\cC$, that is, no nonidentity matrix in $\l\bA\r$ fixes  $\cC$ point-wise. 

\begin{lemma}\label{irreducible}
Let $\cC$ be a $\l \bA\r$-irreducible linear $[n,k]_q$ code. Then either $k=2$, or $k=1$, $q$ is even and $\cC$ is  spanned by the vector $(1,1,\ldots,1)$.
If further $\cC$ is faithful, then $\l \bA\r$ is regular on the nonzero codewords;
in particular, $\cC$ is an equidistant $[n,2]_q$ code of weight $q$.
\end{lemma}
\proof
Assume that $\bA$ induces an invertible linear transformation of order $m$ on  $\cC$. Then $m>1$, and $m$ is a divisor of $q^2-1$. Now $k$ is the smallest positive integer such that $q^k-1\equiv 0\,(\mod m)$, refer to \cite[p.165, II.3.10]{Huppert}.  Thus $k\leqslant 2$.

\vskip 5pt

Suppose that $k=1$. Then $m$ is a divisor of $q-1$, and the
kernel of $\l \bA\r$ acting  on $\cC$ contains the  unique subgroup $\l \bA^{q-1}\r$
of order $q+1$.
Thus $\bu\bA^{q-1}=\bu$ for all $\bu \in \cC$. If $q$ is odd, then $(\bA^{q-1})^{q+1\over 2}=(\bA^n)^{q-1\over 2}=(\eta\lambda^2)^{q-1\over 2}\bI_n=-\bI_n$, yielding $\bu=\bu(\bA^{q-1})^{q+1\over 2}=-\bu$,
which is impossible. Therefore, $q$ is a even. For an arbitrary codeword $\bu=(u_1,u_2,\ldots,u_n)\in \cC$, calculation shows that
\[(u_1,u_2,\ldots,u_n)\bA^{q-1}=(u_3, u_4, u_5,\ldots, u_n,u_1,u_2).\]
Since $\bu\bA^{q-1}=\bu$, we have $u_1=u_2=\cdots=u_n$. Then $\cC$ is spanned by the vector $(1,1,\ldots,1)$, and the first part of this lemma follows.

\vskip 5pt

Now let $\cC$ be faithful. Then $k=2$. Noticing the Singleton bound, we may choose a nonzero word $\bw_1$ with $\wt(\bw_1)\leq n-1$. Let
$\bw_2=\bw_1\bA$. Recalling that $\cC_{\bw_1}$ is $1$-dimensional, it is not $\l \bA\r$-invariant, and thus $\cC_{\bw_1}\ne \cC_{\bw_2}$. In particular, $\cC=\cC_{\bw_1}\oplus \cC_{\bw_2}$.
Assume that $\bA^i$ fixes $\bw_1$ for some $i$. Then $\bA^i$ also fixes $\bw_2$, and so $\bA^i$ fixes $\cC$ point-wise. This implies that $\bA^i=\bI_n$. Then $\l \bA\r$ is regular on $\cC\setminus\{{\bf 0}\}$, and the lemma follows from Lemma \ref{eqid}
\qed

\vskip 5pt

\begin{theorem}\label{EDC}
Assume that $n=q+1=2^sr^t$ for some odd prime $r$ and integers $s,t\geqslant 0$.
Then there exists a faithful $\l \bA\r$-irreducible liner $[n,2]_q$ code. If   $q$ is a Merdenne prime then $\FF_q^n$ is a direct sum of faithful $\l \bA\r$-irreducible linear $[n,2]_q$ codes.
\end{theorem}
\proof
%Let $\cC$ be a $\l \bA\r$-irreducible $[n,k]_q$ code.
Appealing to Maschke's Theorem, refer to \cite[p.123, I.17.7]{Huppert},  we write
\[\FF_q^n=\oplus_{i=1}^m\cC_i,\]
where   $\cC_i$ are   $\l \bA\r$-irreducible $[n,k_i]_q$ codes.
By Lemma \ref{irreducible}, we  assume that $k_1=\cdots=k_{m-1}=2$, and either $k_m=2$ or $q$ is even and $\cC_m$ is spanned by   $(1,1,\ldots,1)$.

\vskip 5pt

Let $K_i$ be the kernel of $\l \bA\r$ acting on $\cC_i$, where
$1\leqslant i\leqslant m$. Recalling that $\bA^n=\eta\lambda^2\bI_n$, we know that $\l\bA^n\r$ is semiregular on the set of nonzero codewords of every   $\cC_i$, and thus  $K_i\cap\l\bA^n\r=1$. Then $|K_i|$ is a  divisor of $q+1$. Now it suffices to show that $|K_i|=1$ for some $i$, and if $q$ is a Merdenne prime then $|K_i|=1$ for all $i$.

\vskip 5pt

Assume first $q$ is even. Then $n=r^t$, and $\l \bA\r$ contains a unique subgroup of order $r$. It follows that either $|K_i|=1$ for some $i$, or all $K_i$ contains a common subgroup of order $r$.
The latter case implies that $\l \bA\r$ is unfaithful on $\FF_q^n$, which is impossible.

\vskip 5pt

Now let $q$ be odd. Then $\l \bA^n\r$ has even order $q-1$. Recalling that $K_i\cap\l\bA^n\r=1$ for all $i$, since $\l \bA\r$ has a unique involution, it follows that every $|K_i|$ is an odd  divisor of $q+1$. Thus, since $\l \bA\r$ is faithful on $\FF_q^n$, we have $|K_i|=1$ for some $i$. If further $q$ is a Merdenne prime, then $|K_i|=1$ for all $i$. This completes the proof.
\qed

\vskip 20pt

\section{A construction of graphs with non-diagonal PA type}\label{sect=exam}

For  a finite group $G$  and    $H\leqslant G$, denote by $[G:H]$ the set   of right cosets of $H$ in $G$. Assume that $H$ is core-free in $G$, that is, $\cap_{g\in G}H^g=1$.
 Then we have a faithful and transitive action of $G$  on $[G:H]$  by right multiplication, and thus we  identify $G$ with a transitive permutation group on $[G:H]$. For a $2$-element $g\in G\setminus H$ with $g^2\in H$,  the coset graph $\Cos(G,H,g)$ is defined as the graph with vertex set $[G:H]$ such that  $Hx$ and $Hy$ are adjacent if and only if $yx^{-1}\in HgH$. It is well-known that   $\Cos(G,H,g)$ is   $G$-arc-transitive and of valency  $|H:(H\cap H^g)|$, and that up to isomorphism every arc-transitive graph is constructed in this way.
As a graph automorphism, the element $g$ maps the vertex $H$ to one of its neighbors, it follows that $\Cos(G,H,g)$  is connected if and only if $G=\l H,g\r$, refer to \cite[p.118, 17B]{Biggs}.

 \vskip 5pt

In the following, for some prime power $q$, we will construct a quasiprimitive group $G$ of  (non-diagonal) PA type with a point stabilizer $H$ isomorphic to the affine group $\AGL_1(q^2)$, and then produce a connected coset graph $\Cos(G,H,g)$ of valency $q^2$. If this is so then, noting that $H$ acts $2$-transitively on $[H:(H\cap H^g)]$ by right multiplication, $\Cos(G,H,g)$ is $(G,2)$-arc-transitive by \cite[Theorem 2.1]{FP1}; of course,  such a graph satisfies Theorem \ref{th-1} (2).

\vskip 10pt

For the rest of this section, we always assume that
%\begin{hypothesis}\label{hypo}
\begin{itemize}
  \item[(C1)] $q=p^f$ for some prime $p$ and integer $f\geqslant 1$, and $n:=q+1=2^sr^t>3$, where  $t\geqslant 0$, $r$ is an odd prime, and either $s\geqslant 2$ or $q$ is even;
  %,  and $n:=q+1$ is either a prime power or the product
  %of two prime powers;
  \item[(C2)] $X$ is an almost simple group with socle $T$,  $|X:T|\leqslant 2$ and
  $X$ has a subgroup $R$ isomorphic to $\AGL_1(q)$, write $R=F{:}(\l b\r\times \l c\r)$, where $F\cong \ZZ_p^f$,
  $b$ has order ${q-1\over (2,q-1)}$ and $c$ has order $(2,q-1)$;
\item[(C3)] $\tau=(1,2,\ldots,n)\in \Sy_n$, and $W=X\wr \l\tau\r$, the wreath product of $X$ by $\l\tau\r$, where
    \[(x_1,x_2,\ldots,x_n)^\tau=(x_n,x_1,x_2,\ldots,x_{n-1})  \,\mbox{ for } x_i\in X, \, 1\leqslant i\leqslant n;\]
%$X^n=\{(x_1,x_2,\ldots,x_n)\mid x_j\in X\}$ and
%$T^n=\{(t_1,t_2,\ldots,t_n)\mid t_j\in T\}$,
%the (external) direct products of $n$ copies of $X$
%and $T$, respectively;
%\item[(iv)] $T_i=\{(t_1,t_2,\ldots,t_n)\mid t_j\in T, t_j=1 %\mbox{ if } j\ne i, 1\leqslant j\leqslant n\}$, where %$1\leqslant i\leqslant n$;
\item[(C4)] $\pi_i: (x_1,x_2,\ldots,x_n)\mapsto x_i$, $1\leqslant i\leqslant n$, are the projections of $X^n$ onto $X$.
\end{itemize}
%\end{hypothesis}

 %Write the elements of $W$ in the form of ${\bf x}\tau^i$
  %with ${\bf x}\in X^n$, and view $\tau$ as an element of $W$.

\vskip 10pt

The next lemma follows easily from (C1) and (C2).

\begin{lemma}\label{F<T}
$R\cap T= F{:}\l b,c^{|X:T|}\r$.
\end{lemma}
\proof
Note that $F$ is the unique minimal normal subgroup of $R$.   Since $T\unlhd X$, we have $F\cap T\unlhd R$, yielding $F\cap T=1$ or $F\leqslant T$. If $F\cap T=1$ then $|X|$ is divisible
by $|F||T|=q|T|$, yielding $|X:T|\geqslant q\geqslant 3$, a contradiction. Thus $F\leqslant T$. Since $|X:T|\leqslant 2$, we have $b\in T$.
Then \[R\cap T=F{:}(\l b,c\r\cap T)=F\l b\r(\l c\r\cap T)= F{:}\l b,c^{|X:T|}\r,\] as desired. This completes the proof.
\qed

\vskip 5pt

For  $Y\leqslant X$, we always deal with   the direct product $Y^n$ of $n$ copies $Y$  as a subgroup  of $W$. Also,
$\l \tau\r$ is viewed as a subgroup of $W$,
so that $W=X^n{:}\l\tau\r$. Sometimes, we
 use boldface type for the elements in $X^n$.
 %For example, ${\bf 1}=(1,1,\ldots, 1)$,
 %the identity of $X^n$.
 Pick three  elements in   $R^n$   as follows:
 \[
 \bb=(b,b,\ldots,b),\,\bc=(c,c,\ldots,c),\,
 \bd_0=(bc,b,bc,\ldots,bc).
 \]
 Then $\bb$, $\bc$ and $\bd_0$ have order ${q-1\over (2,q-1)}$, $(2,q-1)$ and $q-1$, respectively.
 Let \[\theta=\bd_0\tau.\]
 Then\[
 \theta^n= (b^2c,b^2c,b^2c,\ldots,b^2c)=\bb^2\bc.\]
It follows that $\theta$ has order $n(q-1)=q^2-1$, and $\l \theta^n\r=\l \bb^2\r\times \l\bc\r$.

\vskip 5pt

It is easy to check that $\bC_{\l\theta\r}(F^n)=1$,
$F^n\cap \l \theta\r=1$ and $F^n$ is normalized by $\theta$.
Viewing $F^n$ as the $n$-dimensional vector space  $\FF_q^n$, by Lemma \ref{eqid} and Theorem \ref{EDC}, we have the following lemma.

\begin{lemma}\label{tech}
$F^n{:}\l \theta\r$ has a minimal normal subgroup $E$ such that
\begin{itemize}
  \item[(1)] $E\cong \ZZ_p^{2f}$;
  \item[(2)] $\l \theta\r$ acts transitively on $E\setminus\{1\}$ by conjugation, in particular,  $E{:}\l \theta\r\cong \AGL_1(q^2)$;
  \item[(3)] $\pi_i(E)=F$ and $\ker(\pi_i)\cap E\ne \ker(\pi_j)\cap E$, where $1\leqslant i<j\leqslant n$.
\end{itemize}
\end{lemma}

\vskip 5pt

Using Lemma \ref{tech}, we can easily construct a quasiprimitive permutation group   of PA type, which is described as in the  following result.

\begin{theorem}\label{qusi-p-g}
%Assume that $n=q+1=2^ur^v$ for some integers $u,v\geqslant 0$.
Let $G=T^n\l\theta\r$, and let $E$ be a minimal normal subgroup of $F^n{:}\l \theta\r$ satisfying {\rm (1)-(3)} of Lemma {\rm \ref{tech}}.
Let $H=E{:}\l \theta\r$. Then $G$ is a quasiprimitive group on $[G:H]$ of {\rm (}non-diagonal{\rm )} {\rm PA} type, where $T^n\cap H$ is a subdirect product of $(R\cap T)^n$.
\end{theorem}
\proof
First, it is easily shown that $\bC_{\l \theta\r}(T^n)=1$, and $\l \theta\r$ normalizes $T^n$ and acts transitively by conjugation on the set of simple direct factors of $T^n$.
This implies that $G$ is a group and has a unique   minimal normal subgroup $T^n$, and hence $H$ is core-free in  $G$. Thus it suffices to show that $\pi_i(T^n\cap H)=R\cap T$ for $1\leqslant i\leqslant n$.

\vskip 5pt

Calculation shows that $\theta^m\in T^n$ if and only if $m$ is divisible by $n|X:T|$. It follows that $T^n\cap \l\theta\r=\l \theta^{n|X:T|}\r=\l \bb^{2|X:T|}, \bc^{|X:T|}\r$. Since either $q\equiv -1\,(\mod 4)$ or $q$ is even, $\bb$ has odd order ${q-1\over (2,q-1)}$. Noting that
$2|X:T|$ is a divisor of $4$, we have $\l \bb^{2|X:T|}\r=\l \bb\r$. Then $T^n\cap \l\theta\r=\l \bb, \bc^{|X:T|}\r$. Now \[T^n\cap H=T^n\cap (E{:}\l \theta\r)=E{:}(T^n\cap \l \theta\r)=E{:}(\l \bb, \bc^{|X:T|}\r). \]
By Lemmas \ref{F<T} and \ref{tech},  we have
\[R\cap T=F(\l  b,c^{|X:T|}\r)=\pi_i(E)\pi_i(\l \bb, \bc^{|X:T|}\r)=\pi_i(T^n\cap H),\] as desired.
This completes the proof.
\qed

\vskip 10pt

Now we are ready to give a construction for graphs of non-diagonal PA type.
\begin{theorem}\label{const}
Let $G$ and $H$ be  as in Theorem {\rm \ref{qusi-p-g}}.  Suppose that $\N_X(\l b, c\r)$ contains an involution $o$ of $T$ such that   $X=\l F,b,c,o\r$. Let ${\bf o}=(o,o,\ldots, o)$ and
$\Ga(X)=\Cos(G,H,{\bf o})$. Then $\Ga(X)$ is   connected, $(G,2)$-arc-transitive and of valency $q^2$.
\end{theorem}
\proof
We first show that $\Ga(X)$ is  $(G,2)$-arc-transitive. Noting that $H\cong \AGL_1(q^2)$, if  $H\cap H^\bo$ has order $q^2-1$ then $\Ga(X)$ have valency $q^2$, which yields the $2$-arc-transitivity of $G$ on the  graph $\Ga(X)$. Thus it suffices to confirm that $|H\cap H^\bo|=q-1$ and $G=\l H,\bo\r$.

\vskip 5pt

By the choice of $o$, we know that
$o$ centralizes $c$ and normalizes $\l b\r$. Let
%${\bf b}=(b,b,\ldots, b)$ and
${\bf c}_0=(c,1,c,\ldots, c)$.
Then $\bo$ centralizes
${\bf c}_0$ and normalizes $\l {\bf b}\r$. Clearly, $\bo$ centralizes $\tau$. Then $\bo$ centralizes ${\bf c}_0\tau$.
Noting that $\l \theta\r=\l \bd_0\tau\r=\l {\bf b}\r\times \l {\bf c}_0\tau\r$, it follows that $\l \theta\r^{\bo}=\l \theta\r$, and so $\l \theta\r\leqslant H\cap H^\bo$. Suppose that $|H\cap H^\bo|>q-1$.
Since $\l \theta\r$ is maximal in $H$, we have $H\cap H^\bo=H$,
which yields that $E$ is normalized by $\bf o$.
Then $\pi_1(E)=F$ is normalized by $o$. Since $\l F,b,c,o\r=X$, we have $F\unlhd X$, which is impossible. Thus $|H\cap H^\bo|=q-1$, as desired.

\vskip 5pt

By the choice of $(X,T,o)$, we have $T=\l F,b,c^{|X:T|},o\r$. Recalling that $\theta^n={\bf b}^2{\bf c}$, since $\bf b$ has odd order, we have $\l \theta^n\r=\l {\bf b}\r\times \l {\bf c}\r$. By Lemma \ref{tech}, $\pi_i(E)=F$ for all $i$. We have $\pi_i(T^n\cap \l H,\bo\r)\geqslant \l F,b,c^{|X:T|},o\r=T$, yielding $\pi_i(T^n\cap \l H,\bo\r)=T$, where $1\leqslant i\leqslant n$. Let $K_i=\ker(\pi_i)\cap T^n$. Then
\[(T^n\cap \l H,\bo\r)/K_i\cong T,\,\,1\leqslant i\leqslant n.\]
Again by Lemma \ref{tech}, $\ker(\pi_1)\cap E$,  $\ker(\pi_2)\cap E,\ldots, \ker(\pi_n)\cap E$ are distinct. Then $K_1,\ldots,K_n$ are distinct normal subgroups of $T^n\cap \l H,\bo\r$. It follows that $T^n\cap \l H,\bo\r\cong T^n$, refer to \cite[p.113, Lemma 4.3A]{Dixon}. Then $T^n\cap \l H,\bo\r=T^n$, and so
$\l H,\bo\r\geqslant \l T^n,\theta\r=G$. Thus $G=\l H,\bo\r$ as desired. This completes the proof.
\qed

\vskip 10pt

The following example collects some almost simple groups, which support Theorem \ref{const}. Thus there do exist $2$-arc-transitive graphs which satisfy (2) of Theorem \ref{th-1}.

\begin{example}\label{exam}
(1) Let $X=\Sy_p$ and $T=\A_p$, where $7\leqslant p\equiv -1\,(\mod 4)$, and $p+1$ has at most two distinct prime divisors.
Then $\Sy_p$ has a maximal subgroup $F{:}\l a\r$ isomorphic $\AGL_1(p)$, refer to \cite{maxAn}, where $F\cong \ZZ_p$, and $a$ is a $(p-1)$-cycle.
Let $b=a^2$ and $c=a^{p-1\over 2}$. Then $F\l b\r$ is a maximal subgroup of $\A_p$, and $c$ is a product of ${p-1\over 2}$ disjoint transpositions.
It is easy to see that $\Sy_p$ contains an element $d$, which is a product of ${p-1\over 2}$ disjoint transpositions and inverses $a$ by conjugation. Clearly, $cd=dc$. Let $o=cd$. We have $oc=co$, $o\in \A_n$, and $\l F, b, o\r=\A_n$. Thus, by Theorem \ref{const}, we get
a connected $2$-arc-transitive graph $\Ga(X)$ of valency $p^2$.

\vskip 5pt

(2) Let $X=\PGL_2(q)$ and $T=\PSL_2(q)$, where either $q\geqslant 4$ is even or $7\leqslant q\equiv -1\,(\mod 4)$, and $q+1$ has at most two distinct prime divisors.
Note that all subgroups of $X$ and $T$ are explicitly known, refer to \cite{3-D} and \cite[p.213, II.8.27]{Huppert}, respectively.
In particular, $X$ has a maximal subgroup $F{:}\l a\r$ isomorphic $\AGL_1(q)$, where $|F|=q$, and $a$ has order $q-1$. Let $b=a^{(2,q-1)}$ and $c=a^{q-1\over (2,q-1)}$. Then $F\l b\r$ is a maximal subgroup of $T$.
Let $N=\N_X(\l a\r)$. Then $N$ is a dihedral group of order $2(q-1)$, and  $N\cap T$ is a dihedral group of   order $2(q-1)\over (2,q-1)$. Pick an involution $o$ in  $N\cap T$. Then $oc=co$ and $T=\l F,b,o\r$. By Theorem \ref{const}, we get
a connected $2$-arc-transitive graph $\Ga(X)$ of valency $q^2$.
\qed
\end{example}

\vskip 5pt

We  end this section by an example, which gives some graphs satisfying  Theorem \ref{th-1}~(ii).
\begin{example}\label{2^6}
Let $\PSL_2(8)=T<X=T.3\cong \Ree(3)$, and let $F$ be a Sylow $2$-subgroup of $T$. By the Atlas \cite{Atlas}, we have $\N_T(F)\cong \ZZ_2^3{:}\ZZ_7$ and $\N_X(F)\cong \ZZ_2^3{:}(\ZZ_7{:}\ZZ_3)$. Pick an element $b$ of order $3$ in $\N_X(F)$. Let $\tau$ be the $21$-cycle $(1,2,\ldots,21)$ in $\Sy_{21}$.

\vskip 5pt

It is easily shown that the wreath product $X\wr \l\tau\r$ has a normal subgroup $G=T^{21}{:}\l \theta\r$, where $\theta=(b,1,b,\ldots,d)\tau$ has order $63$.
Let $M=T^{21}$. Then $M$ is the unique minimal normal subgroup of $G$. Note that $F^{21}$ is a $\l\theta\r$-invariant subgroup of $M$.
Considering the conjugation of $\l\theta\r$ on $F^{21}$, calculation with GAP \cite{GAP} shows that
\begin{itemize}
  \item[(1)] $F^{21}$ has exactly $13$ minimal $\l\theta\r$-invariant subgroups: one of them has order $2$, one of them has order $2^2$, two of them have order $2^3$, and the other ones have order $2^6$; in fact, $F^{21}$ is the direct product of these $13$ subgroups;
  \item[(2)]  among those $9$ subgroups of order $2^6$ in (1), there are exactly $6$  subgroups such that $\l\theta\r$ acts regularly on the nonidentity elements, that is, each of these $6$ subgroups together with $\theta$ generates a group isomorphic to   $\AGL_1(2^6)$.
\end{itemize}

\vskip 5pt

We fix a minimal $\l\theta\r$-invariant subgroup $E$ of $F^{21}$ with $E\l\theta\r\cong \AGL_1(2^6)$, and let $H=E\l\theta\r$.
Then $M\cap H=E\cong \ZZ_2^6$, and $G$ is a quasiprimitive group of (non-diagonal) PA type on $[G:H]$.
Consider the normalizer of $\l\theta\r$ in $G$. We have
$\N_G(\l\theta\r)=\N_M(\l\theta\r)\l\theta\r$.
Again  confirmed by GAP \cite{GAP}, we conclude that $\N_M(\l\theta\r)\cong \Sy_3$, $\N_G(\l\theta\r)=\N_M(\l\theta\r)\times \l\theta\r$, and there is a unique $2$-element $g\in \N_M(\l\theta\r)$ (up to the double coset $HgH$) such that $G=\l H,g\r$. Thus we have a connected
$(G,2)$-arc-transitive graph $\Cos(G,H,g)$ of valency $2^6$ and order
$2^{57}\cdot 3^{42}\cdot 7^{21}$, where $M$ acts regularly on the arc set of this graph.

\vskip 5pt

Note, there are $6$ choices for the group $E$, and so we may obtain $6$ graphs. However, we do not know  whether there are isomorphic ones among these graphs.
\qed

\end{example}

%%%%%%%%%%%%%%%%%%%%%%%%%%%

\vskip 20pt

\section{A construction of bipartite graphs with diagonal PA type}\label{biparte-exam}

We say a graph is  a  standard double cover
if it is  isomorphic the standard double cover of some graph.
This section aims to construct some  $2$-arc-transitive bipartite graphs with  diagonal PA type, which are not  standard double covers.

\begin{lemma}\label{sdc}
Let $\Ga=(V,E)$ is a connected bipartite graph, $G\leqslant \Aut(\Ga)$. Let $G^*$ be the bipartition preserving subgroup of $G$. Assume that $G$ is transitive on $V$. If $\Ga$ is a  standard double cover, then $\{G_\a\mid \a\in V\}$ is a conjugacy class of subgroups in $G^*$.
\end{lemma}
\proof
Clearly, $G_\a\leqslant G^*$ for all $\a\in V$.
Let $U$ and $W$ be the $G^*$-orbits on $V$. Then
$\{G_\a\mid \a\in U\}$ and $\{G_\b\mid \b\in W\}$ are conjugacy classes of subgroups in $G^*$.
Assume that $\Ga$ is a  standard double cover. 
Then $\Aut(\Ga)$ has an involution $\iota$ which centralizing $G^*$ and interchanges $U$ and $W$. Let $\a\in U$ and $\b=\a^\iota$. We have $\b\in W$. Replacing $G$ by $G^*\times\l\iota\r$ if necessary, we have $G_\b=G_{\a^\iota}=G_\a^\iota=G_\a$. It follows that $\{G_\a\mid \a\in U\}=\{G_\b\mid \b\in W\}$, and the lemma follows.
\qed

\vskip 10pt

From now on, let $p\geqslant 5$ be a prime, and let $\tau=(1,2,\ldots,p-1)\in \Sy_{p-1}$. Let $X=\PGL(2,p)$ or $\Sy_p$ with socle $T$. We will define a subgroup $G$ of
the wreath product $W=X\wr\l\tau\r$, and  construct connected $(G,2)$-arc-transitive bipartite graphs.

\vskip 5pt

Note that $X$ has a   subgroup $R$ isomorphic to $\AGL_1(p)$, and $T\cap R\cong \ZZ_p{:}\ZZ_{p-1\over 2}$.
Choose    $a,b\in R$ with order $p$  and  $p-1$, respectively.
 Then $R=\l a\r{:}\l b\r$.
It is easily shown that $b$ is  contained in a dihedral subgroup $D$ of $X$ with
order $2(p-1)$, which has the center $\l b^{p-1\over 2}\r$ and intersects with $T$ at a dihedral group   of order $p-1$.
Thus both $T$ and $X\setminus T$ contain involutions which
inverse $b$ and centralize $b^{p-1\over 2}$.
Choose  an involution
$c\in X$ with $b^{c}=b^{-1}$ and $cb^{p-1\over 2}\not\in T$.
% such that  and $xc\not\in T$.
We have $D=\l b,{c}\r$, and $X=\l a,b,{c}\r$.

\vskip 5pt

Pick three elements in $W$ as follows:
\[
 \ba=(a,a,\ldots,a),\,
 \bb=(b,b,\ldots,b),\,
\bo=({c},b{c},b^2{c},\ldots,b^{p-2}{c}). 
\]
Clearly, $\tau$ centralizes both $\ba$ and $\bb$, and all coordinates of $\bo$ are distinct.
%, and  $\l \bb,\tau,\bo\r$ normalizes $T^{p-1}$. 
In addition, \[\bo^\tau=\bb^{-1}\bo,\,\bb^\bo=\bb^{-1},\, \tau^\bo=\bb^{-1}\tau,\,\l \ba,\bb,\tau\r=\l\ba\r{:}\l \bb\r\times \l\tau\r,\,  \l\ba,\bb\r\cap T^{p-1}=\l\ba, \bb^2\r.\]
Let
\[G^*=T^{p-1}\l \bb,\tau\r.\]
Suppose that  $\bo \in G^*$. We have  $\bo=(t_1,t_2,\ldots,t_{p-1})\bb^i$ for some $i$ and $t_1,t_2,\ldots,t_{p-1}\in T$. Then $(t_1,t_2,\ldots,t_{p-1})=\bo\bb^{-i}=(b^i{c},b^{i+1}{c},\ldots, b^{p-2+i}{c})$. It follows that $b=b^{i+1}{c}b^i{c}=t_2t_1\in T$, a contradiction. Therefore, $\bo\not\in G^*$.

\vskip 5pt

Let
\[G=G^*{:}\l \bo\r, \, H=\l \ba,\bb,\tau\r.\]
%Then  \[|G|=2|G^*|=2(p-1)|T|^{p-1},\,\, \AGL_1(p)\times %\ZZ_{p-1}\cong \l \ba,\bb\r\times \l\tau\r=H<G^*.\]
Then $\AGL_1(p)\times \ZZ_{p-1}\cong  H<G^*$, and it is easily shown that $T^{p-1}$ is the unique minimal normal subgroup of $G^*$  and  $G$. Thus we have the following lemma.
\begin{lemma}\label{diag-pa}
The group $G$ acts faithfully on $[G:H]$ by right multiplication,
 $G^*$  have two orbits on $[G:H]$, and $G^*$ is a quasiprimitive
 group with diagonal PA type on each of its orbits, $T^{p-1}\cap H=\l\ba,\bb^2\r$ is a diagonal subgroup of $(T\cap R)^{p-1}$.
\end{lemma}

\vskip 5pt

\begin{theorem}
Let $G$, $H$ and $\bo$ be as above, and let  $\Ga=Cos(G,H,\bo)$. Then   $\Ga$ is a connected $(G,2)$-arc-transitive bipartite graph of valency $p$, and $\Ga$ is not a standard double cover.

\end{theorem}
\proof
Let $K=\l \bb,\tau\r$. Then $|H:K|=p$, and $\bo$ normalizes $K$. Thus $H\cap H^{\bo}\geqslant K$. Suppose that $H\cap H^{\bo}> K$. Then
$H=H^\bo$. Noting that $\l \ba\r$ is   characteristic in $H$, it follows that  $\bo$ normalizes $\l \ba\r$, and so $c$ normalizes $\l a\r$. Then $\l a\r\unlhd \l a,b,c\r=X$, a contradiction. Thus $H\cap H^{\bo}=K$. It is easily shown that $H$ acts $2$-transitively on $[H:K]$ by right multiplication. Then $\Ga$ is  $(G,2)$-arc-transitive and of valency $p$.

\vskip 5pt

We next show that $\Ga$ is connected, that is, $G=\l H,\bo\r$.
Let $G_0=\l\ba,\bb, {\bo}\r$.  Clearly,  $G_0$ is a subgroup
of $X^{p-1}$ and normalized by $\tau$. We have $G_0{:}\l\tau\r=\l\ba,\bb,\tau, {\bo}\r=\l H, {\bo}\r$. Then it suffices to show $T^{p-1}\leqslant G_0$.

\vskip 5pt

For $x\in X$, denote by $\be_{i,x}$ the element of $X^{p-1}$ with the $i$th coordinate $x$ and all other coordinates $1$. Write $X^{p-1}=X_1\times X_2\times\cdots \times X_{p-1}$ and  $T^{p-1}=T_1\times T_2\times\cdots \times T_{p-1}$, where
\[
X_i=\{\be_{i,x}\mid x\in X\},\,
T_i=\{\be_{i,t}\mid t\in T\}, \, 1\leqslant i
\leqslant p-1.
\]
For $1\leqslant i<j\leqslant p-1$,
let $\pi_i$ be the projection of $G_0$  to $X_i$,
and define a group homomorphism:
\[\pi_{ij}: G_0\rightarrow X_i\times X_j,\, \be_{1,x_1}\be_{2,x_2}\cdots
\be_{p-1,x_{p-1}}\mapsto \be_{i,x_i}\be_{j,x_j}.
\]
It is easy to see that \[\pi_i(\ker(\pi_j))\times \pi_j(\ker(\pi_i))\leqslant \pi_{ij}(G_0).\] In addition, \[\pi_i(G_0)=\l \be_{i,a},\be_{i,b}, \be_{i,b^{i-1}c}\r=X_i\cong X,\, 1\leqslant i\leqslant p-1.\]

\vskip 5pt

Suppose that $\ker(\pi_i)=\ker(\pi_j)$ for some $1\leqslant i\leqslant j\leqslant p-1$. Define  $\theta: X_i\rightarrow X_j, \,\pi_i(\bx)\mapsto \pi_j(\bx)$, where $\bx$ runs over the elements of $G_0$.
It is easily shown that $\theta$ is a bijection and preserves the operations of groups. Then $\theta$ is an isomorphism, and
\[\theta: \be_{i,a}\mapsto \be_{j,a},\,\be_{i,b}\mapsto\be_{j,b},\, \be_{i,b^{i-1}c} \mapsto\be_{j,b^{j-1}c}.\]
It follows that $X$ has an automorphism $\sigma$ with
\[\sigma: a\mapsto a,\, b\mapsto b,\, b^{i-1}c\mapsto b^{j-1}c.\]
Note that every automorphism of $X$ is induced by the conjugation of some element in $X$. Then there is $x\in X$ such that
\[a^x= a,\, b^x= b,\, (b^{i-1}c)^x= b^{j-1}c.\]
The only possibility is that $x=1$. Then $b^{i-1}c= b^{j-1}c$, yielding $i=j$. Therefore, $\ker(\pi_i)\ne \ker(\pi_j)$ for $1\leqslant i<j\leqslant p-1$.

\vskip 5pt

Recalling that $G_0$ is normalized by $\tau$, it is easily shown that \[(\ker(\pi_i))^\tau=\ker(\pi_{i^\tau}),\, 1\leqslant i\leqslant p-1.\] In particular, we have $\ker(\pi_i)\ne 1$ for all $i$.
%$1\leqslant i\leqslant p-1$.
Let $1\leqslant i<j\leqslant p-1$. Since $\ker(\phi_i)\unlhd G_0$, we have $\pi_j(\ker(\pi_i))\unlhd X_j$. Then either $T_j\le \pi_j(\ker(\pi_i))$, or $\pi_j(\ker(\pi_i))=1$. The latter case implies that $\ker(\pi_i)=\ker(\pi_j)$, a contradiction.
Thus $T_j\le \pi_j(\ker(\pi_i))$. Similarly, we have $T_i\le \pi_i(\ker(\pi_j))$. Then  \[T_i\times T_j\leqslant \pi_i(\ker(\pi_j))\times \pi_j(\ker(\pi_i))\leqslant \pi_{ij}(G_0).\]
By \cite[p.79, Lemma 4.10]{PS}, we have
$T^{p-1}=T_1\times T_2\times\cdots\times T_{p-1}\leqslant G_0$, as desired.

\vskip 5pt

Now $\Ga$ is a connected $(G,2)$-arc-transitive graph of valency $p$. Note that $H\leqslant G^*$, and $G^*$ has two orbits on $[G:H]$, see Lemma \ref{diag-pa}. Then $\Ga$ is bipartite. 
Suppose that $H$ and $H^\bo$ are conjugate in $G^*$. Since
$G^*=T^{p-1}H$, there is some $\bt\in T^{p-1}$ such that $H^\bt=H^\bo$. Note that  $H$ has center $\l\tau\r$, and   $H^\bo$ has center $\l\tau^\bo\r$. Recalling that $\tau^\bo=\bb^{-1}\tau$, we have $(\tau^i)^\bt=\bb^{-1}\tau$ for some integer $i$. Calculation shows that $(\tau^i)^\bt=\bt'\tau^i$ for some $\bt'\in T^{p-1}$. It follows that $\bb^{-1}=\bt'\in T^{p-1}$, yielding $\bb \in T^{p-1}$. 
Then $\l\ba,\bb^2\r=T^{p-1}\cap H\geqslant \l \ba,\bb\r$, which is impossible as $\bb$ has even order $p-1$. Therefore, $H$ and $H^\bo$ are not conjugate in $G^*$. By Lemma \ref{sdc}, $\Ga$ is not   a standard double cover.
This complete the proof.
\qed

\end{document}